\documentclass[12pt, twoside, a4paper]{article}
\usepackage{amssymb,amsmath,amsthm, epsfig}
\usepackage{graphicx}
\usepackage{float}
\usepackage[english]{babel}
\usepackage{hyperref}
\theoremstyle{plain}
\newtheorem{theorem}{Theorem}[section]
\newtheorem{lemma}[theorem]{Lemma}

\newtheorem{corollary}[theorem]{Corollary}

\theoremstyle{definition}
\newtheorem{definition}[theorem]{Definition}
\newtheorem{example}[theorem]{Example}

\theoremstyle{remark}
\newtheorem{remark}{Remark}

\begin{document}
\pagestyle{myheadings}
\title{Blending attractors of Iterated Function Systems}
\author{Oliveira, Elismar R.}
\date{\today}

\maketitle

\begin{abstract}
In this paper we discuss a new method to blend fractal attractors using the code map for the IFS formed by the Hutchinson--Barnsley operators of a finite family of hyperbolic IFSs.  We introduce a parameter called blending coefficient to measure the similarity between the blended set and each one of the original attractors. We also introduce a discrete approximation algorithm and prove a rigorous error estimation used to approximate  these new objects. Several simulation results are provided illustrating our techniques.
\end{abstract}

\noindent\textbf{2020 Mathematics Subject Classification.}
Primary 28A80; Secondary 37C45, 37F35, 65D18, 68U05.

\smallskip
\noindent\textbf{Keywords.}
Iterated function systems, fractal attractors, Hutchinson--Barnsley operators, fractal blending, code map, self-similar sets, approximation algorithm, error estimation, numerical simulation.

\section{Introduction}
   Why blend IFSs? It is a natural to associate figures with known models, indeed several elements in our day life looks like the combination of other basic shapes. However, unless our IFS  attractors lye  in geometrical spaces such as $\mathbb{R}^{n}$, where operations between sets are known, there is no natural way to combine such sets, unless we use the intrinsic structure of the code space to do so. 
   
   This subject is fairly new and has received a lot of attention, using different approaches, mainly form the computational treatment point of view. One can also find approaches from an analytical point of view trying to classify this objects and connect it with classical applications. The following papers, and references therein, exemplify that.
   \bigskip 
   
   Regarding the second point of view, we highlight the paper \cite{Bro19}, which could be seen as a precursor of ours.  More specifically, if $\mathcal{F}$ and $\mathcal{G}$ are iterated function systems, then any infinite word in the symbols $\mathcal{F}$ and $\mathcal{G}$  induces a limit set (which we denote here \textit{a blend} of the attractors of $\mathcal{F}$ and $\mathcal{G}$) they ask whether this Cantor set (the blend) can also be   realized as the limit set of a single iterated function system $\mathcal{H}$. They prove that    if $\mathcal{F}$, $\mathcal{G}$, and $\mathcal{H}$ consist of $C^{1+\alpha}$ diffeomorphisms, then under    some additional constraints on $\mathcal{F}$ and $\mathcal{G}$ the answer is no.   In their work the problem is motivated by the spectral theory of one-dimensional quasicrystals.
   
   \bigskip 
   
   For the first type of approach, the paper \cite{Bow95}  offers a brief tutorial on Iterated Function Systems (IFS) and illustrates how small parameter changes can produce diverse fractal patterns.  The paper \cite{Mar04} presents a rigorous method for morphing 2D affine IFS attractors, ensuring affine stability and shape similarity throughout the transformation, while introducing mathematical tools for stable and consistent fractal blending.   The paper \cite{Bou25} introduces Random Nonlinear Iterated Function Systems (RNIFS), extending classical IFS to nonlinear and stochastic settings, and demonstrating how such transformations yield stable, complex fractal attractors with enhanced geometric richness compared to traditional models. A comparative study with the classical Sierpi\'nski triangle shows that RNIFS transformations preserve global shape while introducing finer geometric detail and higher structural variability.   
   
   \bigskip

   In this paper we address this same objects in a general hyperbolic framework, not necessarily assuming any differentiable structure, and for any number of IFSs.  The cantor set generated by a sequence of IFS is then called blend of the respective attractors. We show that it preserves some geometric similarity inherited from the original attractors and estimate it through the blend coefficient. We also explore some geometrical and topological features of these blends. Finally, using the tools from \cite{dCOS21} we also introduce  the discrete approximation algorithm for these blends and prove rigorous error estimates for their convergence.
   
   The paper is organized as follows:\\
   In Section~\ref{sec: Contractive IFSs} we recall the theory for contractive (hyperbolic) iterated function systems including the existence and characterization of the attractor, code space and code map.
   
    In Section~\ref{sec: Blending IFSs} we use the previous section framework to define a new IFS in the space $K^*(X)$ and the maps given by the Hutchinson--Barnsley operators of a finite family of contractive IFSs. We then define the blend of IFS attractors as the images of the code map for this new attractor. After, we introduce the notions of blending sequence and blend coefficients to study the similarity between the blend and the original attractors.
    
     In Section~\ref{sec: A discrete approximation algorithm for the blend of attractors} we introduce the discretization of the space $X$ , inspired  by \cite{dCOS21}, via a $\varepsilon$-net and prove a rigorous estimation of the error between the blend and a finite discrete approximation.
     
     Finally,   Section~\ref{sec: Examples and applications} has two parts. In the first part we apply the ideas we develop to blend known fractals as the Maple Leaf and the Sierpi\'nski triangle, and also by adding a third IFS from\cite{Oli17}. Several images from simulation,  and the respective blend coefficients, are provided. In the second part we connect the blends with Canright's envelope (\cite{Can94}) and attractor coverings from \cite{AnMi23}, obtaining some distance estimations. 
    
\section{Contractive IFSs}  \label{sec: Contractive IFSs}
In this section we recall the basics on the classical IFS theory. Let $(X,d)$ be a complete metric space. 

\begin{definition}
	We say that $R=(X, f_{j}, j=1, \ldots, n)$ is a (continuous) Iterated Function System (IFS) if each $f_j: X \to X$ is a (continuous) function.
\end{definition}

\begin{definition}\label{def: IfS contractive}
  We say that $R=(X, f_{j}, j=1, \ldots, n)$ is (Banach) contractive if there exists numbers $\lambda_j, j=1, \ldots, n$ such that
  \begin{enumerate}
  	\item $\operatorname{Lip}(f_{j}) <\lambda_j, j=1, \ldots, n$;
  	\item $\lambda_j \in[0,1), j=1, \ldots, n$.
  \end{enumerate}
  The number $\lambda_{R} :=\max \{\lambda_j \in[0,1), j=1, \ldots, n\}$ is called the contractivity constant of $R$.
\end{definition}
From now on we assume that each IFS is contractive (hyperbolic).

\begin{definition}\label{def: Hut-Bar Op level 1}
	 The function $F_{R}: 2^X \to 2^X$ given by, 
	 \[F_{R}(B)=\bigcup_{j=1,\ldots,n} f_j(B)\]
	 is called  the Hutchinson-Barnsley operator of $R$.
\end{definition}

\begin{definition}\label{def: dist Hausd}
	The set of non-empty compact parts of $(X,d)$ is denoted $K^*(X)$ and the function $d_H: K^*(X) \times K^*(X) \to \mathbb{R}$ given by
	\[d_H(A,B):=\max(d(A,B), d(B,A))\]
	is called the Hausdorff--Pompeiu distance.
\end{definition}

The following result is well known from the literature and its proof is widely spreed across several textbooks:
\begin{theorem} \label{thm:main IFS}
	 The following claim are true:
	 \begin{enumerate}
	 	\item $(K^*(X), d_H)$ is a complete metric space (compact if $(X,d)$ is compact);
	 	\item $F_{R}(K^*(X))\subset K^*(X)$;
	 	\item $\operatorname{Lip}(F_{R}) = \lambda_{R} <1$ ( $F_{R}$ is a Banach contraction);
	 	\item There exists a unique self-similar set $A_{R}\in K^*(X)$ ($F_{R}(A_{R})=F_{R}$);
	 	\item For any set $B_0 \in  K^*(X)$ the sequence $B_0, B_1, B_2, \ldots B_k:=F_{R}^k(B_{0})$ converges to $A_{R}$ w.r.t. $d_H$.
	 \end{enumerate}
\end{theorem}
The items (1)-(3) can be founded in \cite{Hut81}, \cite{BarnDem85} or \cite{Bar88} and (4)-(5) are an immediate consequence of the Banach fixed point theorem and items (1)-(3).

\begin{definition}\label{def: attract level1}
	 The unique set $A_{R}$ given by Theorem~\ref{thm:main IFS} is called the attractor of the IFS $R=(X, f_{j}, j=1, \ldots, n)$.
\end{definition}
The existence of attractor has been proved under much less restrictive hypothesis than contractivity for generalized IFS with finitely many maps, see \cite{Str20}, also for possibly infinite, compact, weakly-hyperbolic IFSs, see \cite{AJS16}. 

\begin{definition}\label{def: code space Omega}
	Consider $\Omega:=\{1, \ldots, n\}^{\mathbb{N}}=\{a=(a_1,a_2,...),\; |\; a_i \in \{1, \ldots, n\}\}$ the (full) shift space and a fixed  number $0<\lambda<1$. We introduce the distance
	\[d_{\lambda}(a,b)=\lambda^{k}, k:=\min \{j | a_j\neq b_j\},\]
	if $a \neq b$, otherwise put $d_{\lambda}(a,b)=0$. This metric induces the product topology making $(\Omega, d_{\lambda})$ a compact metric space.  We consider also the continuous transformation $\sigma: \Omega \to \Omega$ called the one-sided shift
	\[\sigma(a_1,a_2,...) = (a_2,a_3,...).\]
\end{definition}

The following result is easy to check and is also well-known from the classical IFS literature.
\begin{theorem}\label{thm: code map level 1}
	Consider, for any $a=(a_1,a_2,...) \in \Omega$ the sequence $D_0:=X$, $D_1:=f_{a_1}(X)$, $D_2:=f_{a_1}(f_{a_2}(X)$, etc. Then,
	\begin{enumerate}
		\item The sequence $D_k$ is nested and each $D_k \in K^*(X)$;
		\item $\displaystyle \pi_R(a) := \lim_{k \to \infty} D_k = \bigcap_{k \geq 0}  D_k$ is a singleton;
		\item For any point $x \in X$ we obtain $\displaystyle \pi_R(a) := \lim_{k \to \infty} f_{a_1}(\cdots(f_{a_k}(x))$;
		\item For any point $a  \in \Omega$ and $j =1,\ldots,n$ we obtain  $\pi_R(a) = f_{a_1}(\pi_R(\sigma(a)))$ and $f_j(\pi_R(a)) = \pi_R( j*a)$, where $j*a=(j,a_1,a_2,\ldots)$.
	\end{enumerate}
\end{theorem}

\begin{definition}\label{def: code map level 1}
	 The space $(\Omega, d_{\Omega})$ is denoted the code-space for the IFS R and the map $\pi_R: \Omega \to X$ is called the code-map.
\end{definition}

The next result characterizes the attractor as the image of the code-map. We will present a proof to illustrate this application.
\begin{theorem}\label{thm: code map level 1 continuity and attractor}
	 Consider $\lambda:=\lambda_R$.  The code-map $\pi_R: \Omega \to X$ is $1$-Lipschitz continuous and $A_R=\pi_R(\Omega)$.   
\end{theorem}
\begin{proof}
	  If  $a, b \in \Omega$ are such $a\neq  b$ and 
	  \[d_{\lambda_R}(a,b)=\lambda_R^{k}, k:=\min \{j | a_j\neq b_j\},\]
	 then from Theorem~\ref{thm: code map level 1}, we know that $\pi_R(a) = f_{a_1}(\cdots(f_{a_k}(\pi_R(\sigma^k(a))))$ and 
	  $\pi_R(b) = f_{a_1}(\cdots(f_{a_k}(\pi_R(\sigma^k(b))))$. Using the fact that $R$ is $\lambda_R$ contractive one obtain
	  \[d_{\lambda_R}(\pi_R(a),\pi_R(b))\leq \lambda_R^{k} \; d_{\lambda_R}(\pi_R(\sigma^k(a)),\pi_R(\sigma^k(b))) \leq C d_{\lambda_R}(a,b),\]
	  for $C=\operatorname{diam}(\Omega)$. That is,  the code-map $\pi_R: \Omega \to X$ is $1$-Lipschitz continuous.
	  
	  For the second part, we notice that $B:=\pi_R(\Omega) \in K^*(X)$  since $\pi_R$ is continuous and $(\Omega,d_{\lambda_R})$ is compact. Applying the Hutchinson--Barnsley operator $F_R$ we obtain
	  \[F_{R}(B)=\bigcup_{j=1,\ldots,n} f_j(\pi_R(\Omega)) = \pi_R(\Omega)=B\]
	 from Theorem~\ref{thm: code map level 1}.
	 
	 From the uniqueness in Theorem~\ref{thm:main IFS}, we know that  if $F_{R}(B)=B$ and $B \in K^*(X)$ then $B=A_R$.	 
	 
\end{proof}

\section{Blending IFSs}\label{sec: Blending IFSs}

We now consider the scenario where a finite family of contractive  IFSs, acting on the same metric space  $R_{i}=(X, f_{j}^{i}, j=1,\ldots,n_{i} )$, for $i=1,\ldots,N$, is provided. Each one have his own attractor $A_{R_{i}}  \in K^*(X)$.   In order to find a way to blend these compact sets we consider a new hyperspace $Y:= K^*(X)$ where $A_{R_{i}}$ are just points. The Hutchinson-Barnsley operators $F_{R_{i}}: Y \to Y$, for $i=1,\ldots,N$ are contractive by Theorem~\ref{thm:main IFS} which allow us to build a new contractive IFS blending the characteristics of each one of the $R_{i}$'s.

\begin{definition} \label{def: blending IFS} Let $(X,d)$ be a complete metric space and $R_{i}=(X, f_{j}^{i}, j=1,\ldots,n_{i} )$, for $i=1,\ldots,N$ a family of contractive IFSs. We denote $\mathcal{R}:=( K^*(X), F_{R_{i}}, i=1,\ldots,N )$ the blending of the $R_{i}$s. 
\end{definition}

It is immediate, from	Theorem~\ref{thm:main IFS}, that 
$\mathcal{R}$ is a contractive IFS itself on the complete metric space  $(K^*(X) , d_H)$ such that each map $F_{R_{i}}: K^*(X) \to K^*(X)$ is a  Banach contraction with $\operatorname{Lip}(F_{R_{i}}) = \lambda_{R_{i}} <1$.

\begin{definition}\label{def: blending space}
	Consider $ \Omega:=\{1, \ldots, N\}^{\mathbb{N}}=\{\theta=(\theta_1,\theta_2,...),\; |\; \theta_i \in \{1, \ldots, N\}\}$ the code space for $\mathcal{R}$ (do not confuse that with the code space $ \Omega_i:=\{1, \ldots, n_i\}^{\mathbb{N}}$ for each individual IFS $R_i$). We denote  $ \Omega$ the blending space and an element $\theta \in  \Omega$ a blending sequence.  The compact set
	\[\mathcal{A}(\theta) :=\pi_{\mathcal{R}}(\theta) = \lim_{k \to \infty}  F_{R_{\theta_1}}(\cdots( F_{R_{\theta_k}}(Z))  \in K^*(X)\]
	for a fixed $Z \in K^*(X)$, given by Theorem~\ref{thm: code map level 1}, is called the blend  of the attractors  $A_{R_{i}}, i=1,\ldots,N$ by $\theta$.
\end{definition} 
At this point some commentaries are in order. 

\begin{remark}\label{rem: blend 1}
 The blending does not depend on $Z \in K^*(X)$ according to Theorem~\ref{thm: code map level 1}.  Each blending is an element of the attractor of $\mathcal{R}$ forming a subset of $K^*(K^*(X))$:
\[A_{\mathcal{R}}= \pi_{\mathcal{R}}( \Omega).\]
Moreover, from Theorem~\ref{thm: code map level 1 continuity and attractor}, the blend of attractors, $\theta \mapsto \mathcal{A}(\theta) $ is continuous with respect to the blending sequence because $\mathcal{A}(\theta) =\pi_{\mathcal{R}}(\theta)$.   In \cite{Bro19} the set $\mathcal{A}(\theta)$ is called the \emph{dynamically defined Cantor set}(w.r.t. the sequence $\theta$).
\end{remark}

\begin{remark}\label{rem: blend 2}
Despite the extreme complexity of $A_{\mathcal{R}}$ as a compact set of compact sets, we can fully understand it through the analysis of each blend $\mathcal{A}(\theta)$ (see the set of blends in Figure~\ref{fig:blend sierp mapl} and Figure~\ref{fig:blend sierp mapl p2}). This is another application of the blend idea.
\end{remark}

\begin{remark} \label{rem: blend 3}
	The word ``blend" is originated from the interpretation of $\mathcal{A}(\theta)$ as the limit of a sequence $F_{R_{\theta_1}}(\cdots( F_{R_{\theta_k}}(Z))$ blending the Hutchinson--Barnsley of different IFSs according to the ``recipe" $\theta=(\theta_1,\theta_2,...)$. It is evident that, if we do not have any blending, that is, a constant sequence $\theta=(\theta_1,\theta_1,\theta_1,\theta_1,\theta_1,\theta_1,...)$, then
\[\mathcal{A}(\theta) :=\lim_{k \to \infty}  F_{R_{\theta_1}}(\cdots( F_{R_{\theta_1}}(Z))  = \lim_{k \to \infty}  F_{R_{\theta_1}}^k(Z)=A_{R_{\theta_1}} \]
by Theorem~\ref{thm:main IFS}. This shows that the usual IFS fractal operator iteration is just a particular case of blends.
\end{remark}

This remark motivates the interpretation of the blend of attractor as the result of a random iteration of the IFSs  $R_{i}=(X, f_{j}^{i}, j=1,\ldots,n_{i} )$, for $i=1,\ldots,N$. As we will see later, this set inherits some characteristics of each individual attractors $A_{R_{i}}, i=1,\ldots,N$. To quantify that we must take into account the frequency of a given symbol in $\theta$ as well as its position in the sequence.

A way to measure the similarity between the blend of attractors $\mathcal{A}(\theta)$ and each individual $A_{R_{i}}, i=1,\ldots,N$ is estimating the distance to each $A_{R_{i}}$.
\begin{lemma}\label{lem: bound iterated blend}
	 Let $\displaystyle \mathcal{A}(\theta) = \lim_{k \to \infty}  F_{R_{\theta_1}}(\cdots( F_{R_{\theta_k}}(Z)) $
	 for a fixed $Z \in K^*(X)$, be  the blend  of the attractors  $A_{R_{i}}, i=1,\ldots,N$ by $\theta$. Then,
	 \[d_H( F_{R_{\theta_1}}(F_{R_{\theta_2}}(\cdots( F_{R_{\theta_k}}(Z))), A_{R_{i_0}}) \leq \gamma_{k+1} d_H(Z, A_{R_{i_0}}) +\] \[+ \sum_{j=1}^{k} \gamma_j d_H(F_{R_{\theta_j}}(A_{R_{i_0}}), A_{R_{i_0}}),\] 
	 where $\gamma_1:=1$ and $\gamma_j:=\lambda_{R_{\theta_1}} \cdots \lambda_{R_{\theta_{j-1}}}, \; j \geq 2$.
\end{lemma}
\begin{proof}
	
We notice that fixed $W$ and $A_{R_{i_0}}$ we have 
\[d_H( F_{R_{\theta_1}}(W), A_{R_{i_0}}) \leq d_H( F_{R_{\theta_1}}(W), F_{R_{\theta_1}}(A_{R_{i_0}})) +d_H( F_{R_{\theta_1}}(A_{R_{i_0}}), A_{R_{i_0}})\leq\]
\[\lambda _{R_{\theta_1}} d_H(W, A_{R_{i_0}}) + d_H(F_{R_{\theta_1}}(A_{R_{i_0}}), A_{R_{i_0}}).\]
Taking $W:=F_{R_{\theta_2}}(\cdots( F_{R_{\theta_k}}(Z))$  we can repeat  this process, obtaining
\[d_H( F_{R_{\theta_1}}(F_{R_{\theta_2}}(\cdots( F_{R_{\theta_k}}(Z))), A_{R_{i_0}}) \leq \] 
\[\leq\lambda_{R_{\theta_1}} d_H(F_{R_{\theta_2}}(\cdots( F_{R_{\theta_k}}(Z)), A_{R_{i_0}}) + d_H(F_{R_{\theta_1}}(A_{R_{i_0}}), A_{R_{i_0}}) \leq \]
\[\leq \lambda_{R_{\theta_1}} \lambda_{R_{\theta_2}} d_H(F_{R_{\theta_3}}(\cdots( F_{R_{\theta_k}}(Z)), A_{R_{i_0}}) + \lambda_{R_{\theta_1}}d_H(F_{R_{\theta_2}}(A_{R_{i_0}}), A_{R_{i_0}})   + \] \[d_H(F_{R_{\theta_1}}(A_{R_{i_0}}), A_{R_{i_0}}) \leq \]
\[\leq \lambda_{R_{\theta_1}} \lambda_{R_{\theta_2}} \lambda_{R_{\theta_3}} d_H(F_{R_{\theta_4}}(\cdots( F_{R_{\theta_k}}(Z)), A_{R_{i_0}}) + \lambda_{R_{\theta_1}}\lambda_{R_{\theta_2}}d_H(F_{R_{\theta_3}}(A_{R_{i_0}}), A_{R_{i_0}}) \] \[+ \lambda_{R_{\theta_1}}d_H(F_{R_{\theta_2}}(A_{R_{i_0}}), A_{R_{i_0}})   +d_H(F_{R_{\theta_1}}(A_{R_{i_0}}), A_{R_{i_0}}) \leq \]
\[\leq \lambda_{R_{\theta_1}} \cdots \lambda_{R_{\theta_k}} d_H(Z, A_{R_{i_0}}) + \lambda_{R_{\theta_1}} \cdots \lambda_{R_{\theta_{k-1}}} d_H(F_{R_{\theta_k}}(A_{R_{i_0}}), A_{R_{i_0}}) + \ldots\] \[ +\lambda_{R_{\theta_1}}\lambda_{R_{\theta_2}}d_H(F_{R_{\theta_3}}(A_{R_{i_0}}), A_{R_{i_0}})  + \lambda_{R_{\theta_1}}d_H(F_{R_{\theta_2}}(A_{R_{i_0}}), A_{R_{i_0}})   +d_H(F_{R_{\theta_1}}(A_{R_{i_0}}), A_{R_{i_0}}).\]

Thus, 
\[d_H( F_{R_{\theta_1}}(F_{R_{\theta_2}}(\cdots( F_{R_{\theta_k}}(Z))), A_{R_{i_0}}) \leq \lambda_{R_{\theta_1}} \cdots \lambda_{R_{\theta_k}} d_H(Z, A_{R_{i_0}}) + \]\[+\sum_{j=1}^{k} \gamma_j d_H(F_{R_{\theta_j}}(A_{R_{i_0}}), A_{R_{i_0}}),\] 
where $\gamma_1:=1$ and $\gamma_j:=\lambda_{R_{\theta_1}} \cdots \lambda_{R_{\theta_{j-1}}}, \; j \geq 2$.
\end{proof}
Notice that, if we choose $Z=A_{R_{i_0}}$ then $d_H(Z, A_{R_{i_0}})=0$. In any case, 
\[\lambda_{R_{\theta_1}} \cdots \lambda_{R_{\theta_k}} \leq \lambda_{\mathcal{R}}^k \to 0\]
when $k\to \infty$ so the second part will be the most significant for that bound.

\begin{definition}\label{def: blendin coeff }
	Let $\theta \in  \Omega$ a blending sequence and 
	$\mathcal{A}(\theta) $ be  the blend  of the attractors  $A_{R_{i}}, i=1,\ldots,N$ by $\theta$. For each $i \in \{1,\ldots,N\}$ we define the blending coefficient of $\theta$ w.r.t. $i$ as the number
	\[\beta(\theta,i):= \gamma_1 + \sum_{k=2,\; \theta_{k}\neq i}^{\infty}  \gamma_k,\]
	where $\gamma_1:=1$ and $\gamma_k:=\lambda_{R_{\theta_1}} \cdots \lambda_{R_{\theta_{k-1}}}, \; k \geq 2$.
\end{definition}

This number $\beta(\theta,i)$ is well-defined, since $\gamma_k \leq \lambda_{\mathcal{R}}^k, \; k \geq 2$ and $\lambda_{\mathcal{R}}<1$.  We notice that when a symbol $\theta_{k}=i$ appears in the blend sequence we obtain $d_H(F_{R_{\theta_j}}(A_{R_{i}}), A_{R_{i}})= d_H(F_{R_{i}}(A_{R_{i}}), A_{R_{i}})= d_H(A_{R_{i}}, A_{R_{i}})= 0$. Hence,  the additive $\gamma_j$ is eliminated and the bound decreases.  In such times when $\theta_{k}=i$, the iteration is actually exponentially approaching $A_{R_{i}}$. Thus, the coefficient $\beta(\theta,i)$ computes the distancing of $\mathcal{A}(\theta)$ from $A_{R_{i}}$.

\begin{theorem}\label{thm: blendin coeff property}
	Let $\theta \in  \Omega$ a blending sequence and 
	$\mathcal{A}(\theta)$ be  the blend  of the attractors  $A_{R_{i}}, i=1,\ldots,N$ by $\theta$.  The following properties hold:
	\begin{enumerate}
		\item $1 \leq \beta(\theta,i)  \leq \frac{1}{1- \lambda_{\mathcal{R}}}$ for each $ i=1,\ldots,N$, and $\beta(\theta,i)=1$ if, and only if, $\theta = (i, i, i, i, i, i,...)$;
		\item For any $\varepsilon>0$ such that 
		\[d_H(\mathcal{A}(\theta), A_{R_{i_0}})  \leq   \beta(\theta,i_0) \delta_{i_0}\] 
		where $\displaystyle \delta_{i_0}:=\max_{i=1,\ldots,N} d_H(F_{R_{i}}(A_{R_{i_0}}), A_{R_{i_0}})$.
	\end{enumerate}
\end{theorem}
\begin{proof}
	   For the first claim we recall that \[\beta(\theta,i):= \gamma_1 + \sum_{k=2,\; \theta_{k}\neq i}^{\infty}  \gamma_k,\]
	   where $\gamma_1:=1$ and $\gamma_k:=\lambda_{R_{\theta_1}} \cdots \lambda_{R_{\theta_{k-1}}}, \; k \geq 2$.
	   Since $\gamma_k \leq \lambda_{\mathcal{R}}^k, \; k \geq 0$ and $\lambda_{\mathcal{R}}<1$ we get 
	   \[\beta(\theta,i)  \leq   \sum_{k=0}^{\infty} \lambda_{\mathcal{R}}^{k}  = \frac{1}{1- \lambda_{\mathcal{R}}}.\]
	   Since, $\gamma_1=1$, we have from definition that for any sequence but $\theta = (i, i, i, i, i, i,...)$ some positive additive will be computed for $\beta(\theta,i)$ making it strictly bigger than $1$. Hence,  $\beta(\theta,i)=1$ if, and only if, $\theta \sim (i, i, i, i, i, i,...)$.
	   
	   For the second claim we choose any $\varepsilon>0$, and from the limit  
	   \[\mathcal{A}(\theta) = \lim_{k \to \infty}  F_{R_{\theta_1}}(\cdots( F_{R_{\theta_k}}(Z)) \]
	   we select $k$ such that 
	   \[d_H(\mathcal{A}(\theta), F_{R_{\theta_1}}(\cdots( F_{R_{\theta_k}}(Z)) ) < \frac{\varepsilon}{2}.\]
	   
	   Suppose that $k$ is big enough to $\gamma_{k+1} d_H(Z, A_{R_{i_0}}) < \frac{\varepsilon}{2}$ which is always possible because 
	   $\gamma_{k+1} \leq \lambda_{\mathcal{R}}^{k+1} \to 0$ when $k \to \infty$.
	   
	   From Lemma~\ref{lem: bound iterated blend} we obtain
	   \[d_H( F_{R_{\theta_1}}(F_{R_{\theta_2}}(\cdots( F_{R_{\theta_k}}(Z))), A_{R_{i_0}}) \leq \] \[ \leq \gamma_{k+1} d_H(Z, A_{R_{i_0}}) + \sum_{j=1}^{k} \gamma_j d_H(F_{R_{\theta_j}}(A_{R_{i_0}}), A_{R_{i_0}})\leq \] 
	   \[\leq \frac{\varepsilon}{2} +\left(  \sum_{j=1, \; \theta_j \neq i_0}^{k} \gamma_j\right)  \max_{i=1,\ldots,N} d_H(F_{R_{i}}(A_{R_{i_0}}), A_{R_{i_0}})\leq \]
	   \[\leq \frac{\varepsilon}{2} +  \beta(\theta,i_0) \delta_{i_0}.\]
	   Now, we apply the triangular inequality,
	   \[d_H(\mathcal{A}(\theta), A_{R_{i_0}})  \leq \] \[\leq d_H(\mathcal{A}(\theta), F_{R_{\theta_1}}(F_{R_{\theta_2}}(\cdots( F_{R_{\theta_k}}(Z))))  +d_H(F_{R_{\theta_1}}(F_{R_{\theta_2}}(\cdots( F_{R_{\theta_k}}(Z))), A_{R_{i_0}})  \leq\]
	   \[\leq  \frac{\varepsilon}{2} + \frac{\varepsilon}{2} +  \beta(\theta,i_0) \delta_{i_0}= \varepsilon +  \beta(\theta,i_0) \delta_{i_0}.\]
	   Since the number $\varepsilon$ is arbitrary we obtain the desired inequality.
	   
\end{proof}

We notice that the number $\displaystyle \delta_{i_0}=\max_{i=1,\ldots,N} d_H(F_{R_{i}}(A_{R_{i_0}}), A_{R_{i_0}})$ is a measure of how much $A_{R_{i_0}}$ is not self-similar with respect to the remaining IFSs operators $F_{R_{i}}$.

\begin{remark}
	We notice that some obvious questions appears from our discussion. For instance, which other super-structures for a family of IFSs could provide interesting information? Can we fully understand the function $\theta \mapsto \operatorname{Dim}_{H}(\mathcal{A}(\theta)) $? Inspired by \cite{Bro19}, which blends are actually attractors of a single IFS? Are there new blend sets which can not be generated by IFS? Are there nontrivial connections with random IFS or sequential IFS? These questions are, among others, indications that this approach has a large scope and could be very fruitful in the future.
\end{remark}

\section{A discrete approximation algorithm for the blend of attractors}\label{sec: A discrete approximation algorithm for the blend of attractors}
In this section we aim to find a rigorous approximation procedure to represent the blend of attractors as a discrete subset of the original space.

We follow the notation and definitions from \cite[Section 4]{dCOS21}.

\begin{definition} Let $(X,d)$ be a compact metric space.
	\begin{enumerate}
		\item A set $\hat X \subset X$  is called a $\varepsilon$-net if $X \subseteq \hat X^{\varepsilon} $ (that is, for every $x \in X$ there exists $y\in \hat X $ such that $d(x,y)\leq \varepsilon$). The number $\varepsilon$ is called the resolution of the $\varepsilon$-net;
		\item A map $r: X \to \hat X $ is called a $\varepsilon$-projection if $r|_{_{\hat X }}=id$ and $d(x, r(x))\leq \varepsilon, \; \forall x \in X$ (that is, $r(X) \subseteq \hat X^{\varepsilon} $).
		\item Given $B \subset X$ we denote $\widehat{B}:=r(B) \subset \hat X$ the discretization of $B$.
		\item We say that $(\hat X, r)$ is a discretization of $X$, if $\hat X \subset X$  is a $\varepsilon$-net  and $r: X \to \hat X $ is  a $\varepsilon$-projection (w.r.t.  $\hat X$).  The discretization is finite if $\hat X$ is a finite set.
	\end{enumerate}
\end{definition}
\begin{lemma} \label{lem: dH discret lemma}
	 Suppose that $(\hat X, r)$ is a finite discretization of $X$, then $d_H(Z,  \widehat{Z}) \leq \varepsilon$ for any  $Z \in K^*(X)$.
\end{lemma}
\begin{proof}
	Consider $Z \in K^*(X)$. Since the discretization is finite $\widehat{Z} \in K^*(X)$. On one hand, $d(z, r(z))\leq \varepsilon, \; \forall z \in Z$ and $r(z) \in \widehat{Z}$, thus $Z \subseteq \widehat{Z}^\varepsilon$.  Reciprocally, $\widehat{Z}\subseteq Z^\varepsilon$ holds by the same reason. Hence, $d_H(Z,  \widehat{Z}) \leq \varepsilon$.
\end{proof}
\begin{theorem}\label{thm: discretization}
	  For a fixed set  $Z \in K^*(X)$ and $\theta \in \Omega$, consider  $\mathcal{A}(\theta)$ the blend  of the attractors  $A_{R_{i}}, i=1,\ldots,N$ by $\theta$.
	  
	  Suppose that $(\hat X, r)$ is a finite discretization of $X$ and define the following sequence of finite subsets of $\hat X$:
	  \[\begin{cases}
	  	  Y_{k}= \widehat{Z}\\
	  	  Y_{j} = \widehat{  F_{R_{\theta_j}}(\widehat{Y_{j+1}}) } , \; 1 \leq j \leq k-1.
	  \end{cases}\]
	  \begin{enumerate}
	  	\item $d_H(F_{R_{\theta_1}}(\cdots( F_{R_{\theta_k}}(Z)), Y_1) \leq \frac{\varepsilon}{1-\lambda_{\mathcal{R}}} $;
	  	\item $d_H(\mathcal{A}(\theta), Y_1) \leq \lambda_{\mathcal{R}}^k \operatorname{diam}(X) + \frac{\varepsilon}{1-\lambda_{\mathcal{R}}}$.
	  \end{enumerate}
\end{theorem}
\begin{proof}
	  For the first part, we consider  $Y_{k}= \widehat{Z} $, clearly $d_H(Z, Y_{k}) =d_H(Z,  \widehat{Z}) \leq \varepsilon$  from Lemma~\ref{lem: dH discret lemma}.  
	  
	  Now consider the compact set  $Y_{k-1}:=\widehat{F_{R_{\theta_k}}(Y_k)}$, then
	  \[d_H( F_{R_{\theta_k}}(Z), Y_{k-1} ) \leq d_H( F_{R_{\theta_k}}(Z), F_{R_{\theta_k}}(Y_k) )  + d_H( F_{R_{\theta_k}}(Y_k), \widehat{F_{R_{\theta_k}}(Y_k)} ) \leq \]
	  \[\leq \lambda_{R_{\theta_k}} d_H(Z, Y_{k})  + \varepsilon \leq \lambda_{R_{\theta_k}}  \varepsilon  + \varepsilon.\]
	  Inductively, one obtain
	  \[d_H(F_{R_{\theta_1}}(\cdots( F_{R_{\theta_k}}(Z)), Y_1) \leq \lambda_{R_{\theta_1}}\cdots\lambda_{R_{\theta_k}}  \varepsilon  + \ldots +\lambda_{R_{\theta_k}}  \varepsilon  + \varepsilon \leq \frac{\varepsilon}{1-\lambda_{\mathcal{R}}}.\]
	  
	  For the second part we notice that, from Theorem~\ref{thm: code map level 1}, one obtain
	  \[\mathcal{A}(\theta) =\pi_{\mathcal{R}}(\theta) = \lim_{k \to \infty}  F_{R_{\theta_1}}(\cdots( F_{R_{\theta_k}}(Z))  =\]
	  \[F_{R_{\theta_1}}(\cdots( F_{R_{\theta_k}}(W)),\]
	  where $W:=\pi_{\mathcal{R}}(\sigma^k(\theta)) \in \mathcal{A}_{\mathcal{R}}$. Then, using the contractivity of each map we obtain
	  \[d_H(\mathcal{A}(\theta), Y_1) \leq d_H(\mathcal{A}(\theta), F_{R_{\theta_1}}(\cdots( F_{R_{\theta_k}}(Z))) + d_H(F_{R_{\theta_1}}(\cdots( F_{R_{\theta_k}}(Z)), Y_1) \leq \]
	  \[=d_H(F_{R_{\theta_1}}(\cdots( F_{R_{\theta_k}}(W)), F_{R_{\theta_1}}(\cdots( F_{R_{\theta_k}}(Z))) + d_H(F_{R_{\theta_1}}(\cdots( F_{R_{\theta_k}}(Z)), Y_1) \leq\]	  
	  \[\leq  \lambda_{R_{\theta_1}}\cdots\lambda_{R_{\theta_k}} d_H(W, Z) + \frac{\varepsilon}{1-\lambda_{\mathcal{R}}}\leq \lambda_{\mathcal{R}}^k \operatorname{diam}(X) + \frac{\varepsilon}{1-\lambda_{\mathcal{R}}}.\]
	  This concludes our proof.
	  
\end{proof}

\begin{definition}
	The discrete set $B:=Y_1$ appearing in Theorem~\ref{thm: discretization} is called the approximation  of the blend of attractors $A_{R_{i}}, i=1,\ldots,N$ by $\theta$ with resolution $\varepsilon$.  
\end{definition}

We notice that Theorem~\ref{thm: discretization} provides an approximation as good  as the $\varepsilon$ since the term $\lambda_{\mathcal{R}}^k \operatorname{diam}(X)$ can be made arbitrarily small increasing $k$ (that is, taking a longer word of $\theta$). However, the second part  will be at best $\frac{\varepsilon}{1-\lambda_{\mathcal{R}}}$ and can only be decreased by refining the $\varepsilon$-net $\hat X$. In this way, it is natural to compute the absolute error in the approximation  as follows:

The following result generalizes somehow the result in  \cite[Section 4]{dCOS21}, with respect to the generation of the attractor picture,  since we can now apply a different IFS at each iteration.
\begin{corollary} \label{cor: error}  Let  $\delta>0$ a fixed number. Then, there exists $k \in \mathbb{N}$ and $\varepsilon>0$ such that, if  $B=Y_1$ the approximation  of the blend of attractors $A_{R_{i}}, i=1,\ldots,N$ by $\theta$ with resolution $\varepsilon$, then $d_H(\mathcal{A}(\theta), B) \leq \delta$.
\end{corollary}
\begin{proof}
	Theorem~\ref{thm: discretization} provides an approximation 
	\[d_H(\mathcal{A}(\theta), Y_1) \leq \lambda_{\mathcal{R}}^k \operatorname{diam}(X) + \frac{\varepsilon}{1-\lambda_{\mathcal{R}}}.\]
	
	Given $\delta>0$ take a $\varepsilon$-net with $\frac{\varepsilon}{1-\lambda_{\mathcal{R}}} < \frac{\delta}{2}$, that is,  
	$\varepsilon < \frac{\delta}{2(1-\lambda_{\mathcal{R}})}$. Then, choose $k \in \mathbb{N}$ such that $\lambda_{\mathcal{R}}^k \operatorname{diam}(X) < \frac{\delta}{2}$ that is, $k > \frac{\ln\frac{\delta}{2 \operatorname{diam}(X)}}{\ln \lambda_{\mathcal{R}}}$. Then, 
	\[d_H(\mathcal{A}(\theta), B) \leq \delta.\]
\end{proof}

The result in Theorem~\ref{thm: discretization} and Corollary~\ref{cor: error} can be synthesized as an algorithm:
\begin{figure}[H]
	{\tt	{\footnotesize   
			\begin{tabbing}
				aaa\=aaa\=aaa\=aaa\=aaa\=aaa\=aaa\=aaa\= \kill
				\> Input: A list of IFSs $R_{i}=(X, f_{j}^{i}, j=1,\ldots,n_{i} )$, for $i=1,\ldots,N$ \\
				\> Input: A blending sequence $\theta \in \Omega$ and a fixed compact set $Z \in K^*(X)$ \\
				\> Input: The error $\delta>0$ \\
				\> Input: Numbers $k > \frac{\ln\frac{\delta}{2 \operatorname{diam}(X)}}{\ln \lambda_{\mathcal{R}}}$ and $\varepsilon < \frac{\delta}{2(1-\lambda_{\mathcal{R}})}$\\
				\> Output: A discrete approximation $Y_1$ of $\mathcal{A}(\theta) $ with error at most $\delta$ \\
				\>  Algorithm: BlendApprox \\
				\>  	\> $Y_{k}= \widehat{Z}$\\
				\> \> {\bf for } for j from k-1 to 1 do:\\
				\> \> \> $Y_{j} := \widehat{  F_{R_{\theta_j}}(\widehat{Y_{j+1}}) }$\\
				\> 	\>  {\bf end loop} 
	\end{tabbing}}}
	\caption{Approximation of $A(\theta)$ via Algorithm  \textbf{BlendApprox}.}\label{fig:blend approx algorithm} 
\end{figure}

\section{Examples and applications}\label{sec: Examples and applications}
\subsection{Examples}
\begin{example} \label{examp: sierp and maple} 
	 Consider $X=[0,1]^2$ and the IFSs $R_1=(X, f_1, f_2,f_3)$ where
	 \[\begin{cases}
	 	f_1(x,y)=  (0.5 x,  0.5 y )\\
	 	f_2(x,y)=  (0.5  x+ 0.5 ,0.5  y )\\
	 	f_3(x,y)=  (0.5  x+ 0.25 , 0.5  y+0.5 )
	 \end{cases}
	 \]
	 and $R_2=(X, f_1, f_2,f_3,f_4)$ where
	 \[\begin{cases}
	 	f_1(x,y)= (0.8 x1+0.1, 0.8 y1+0.04)\\
	 	f_2(x,y)=  (0.5 x1+0.25, 0.5 y1+0.4 )\\
	 	f_3(x,y)=  ( 0.355 x1-0.355 y1+0.266, 0.355 x1+0.355 y1+0.078)\\
	 	f_4(x,y)= ( 0.355 x1+0.355 y1+0.378, -0.355 x1+0.355 y1+0.434)
	 \end{cases}
	 \]
	 $A_{R_1}$ is a Sierpi\'nski like attractor and  $A_{R_2}$ the Maple leaf attractor (see Figure~\ref{fig:sierp}). We consider the blend of these two attractors. It is well-known that $\lambda_{R_1}=0.5$  and  $\lambda_{R_2}=0.8$ thus $\lambda_{\mathcal{R}}=0.8$.   
	 \begin{figure}[H]
	 	\centering
	 	\includegraphics[width=4cm]{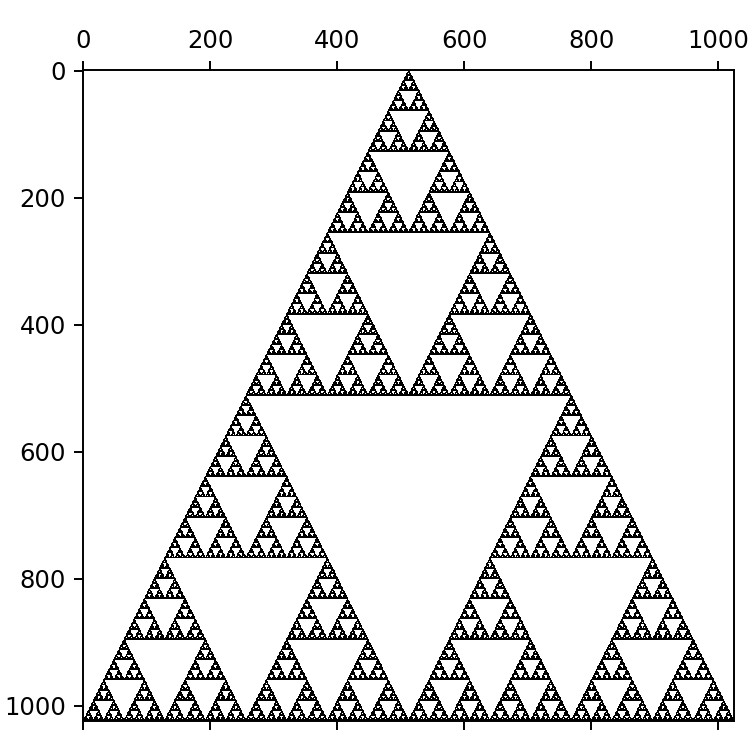} \includegraphics[width=4cm]{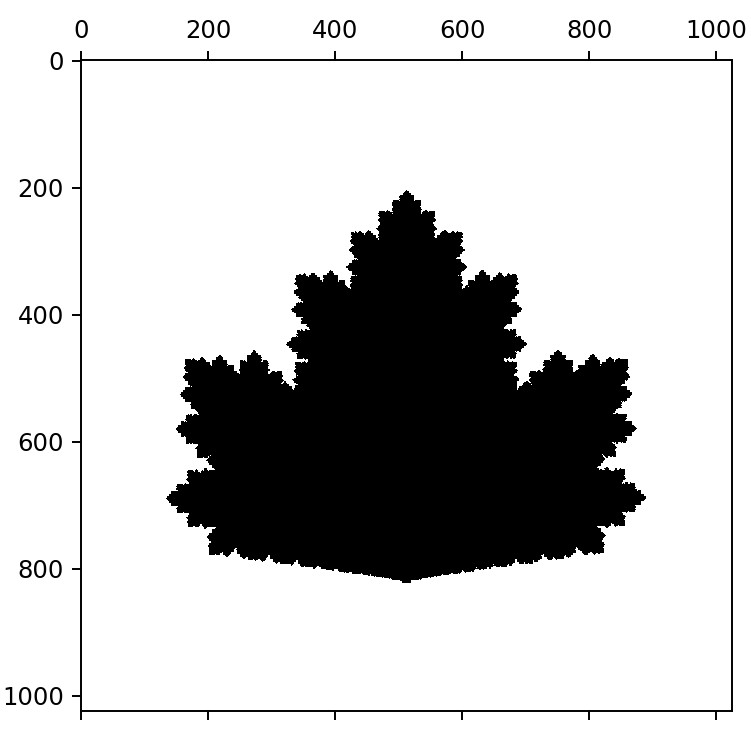}
	 	\caption{Sierpi\'nski(left) and Maple(right).}
	 	\label{fig:sierp}
	 \end{figure}

	 In order to apply the Algorithm~\ref{fig:blend approx algorithm} we consider $M=1024$ and  $\widehat X=\{(x_i,y_j) \in [0,1]^2 \; | \; x_i = 0 + i \frac{1}{M}, \, y_j = 0 + j \frac{1}{M}, \; 0\leq i,j \leq M \}$. Obviously, a $\varepsilon$-net for $\varepsilon=\frac{1}{M \sqrt{2}}$. The $\varepsilon$-projection $r$ is the closest point in  the mesh. We set $k=20$ as  the length of the blending sequence $\theta=(\theta_1,\ldots,\theta_{20})$. 
	 
	 From Theorem~\ref{thm: discretization} one obtain the error as 
	 \[d_H(\mathcal{A}(\theta), Y_1) \leq \lambda_{\mathcal{R}}^k \operatorname{diam}(X) + \frac{\varepsilon}{1-\lambda_{\mathcal{R}}} \approx 0.0198,\]
	 because $\operatorname{diam}(X)=\sqrt{2}$. 
	 
	 We now exhibit some blend discretizations $Y_1$, together with the blend sequence that generate it, and the blend coefficients. 
	 \begin{figure}[H]
	 	\centering
	 	\includegraphics[width=4cm]{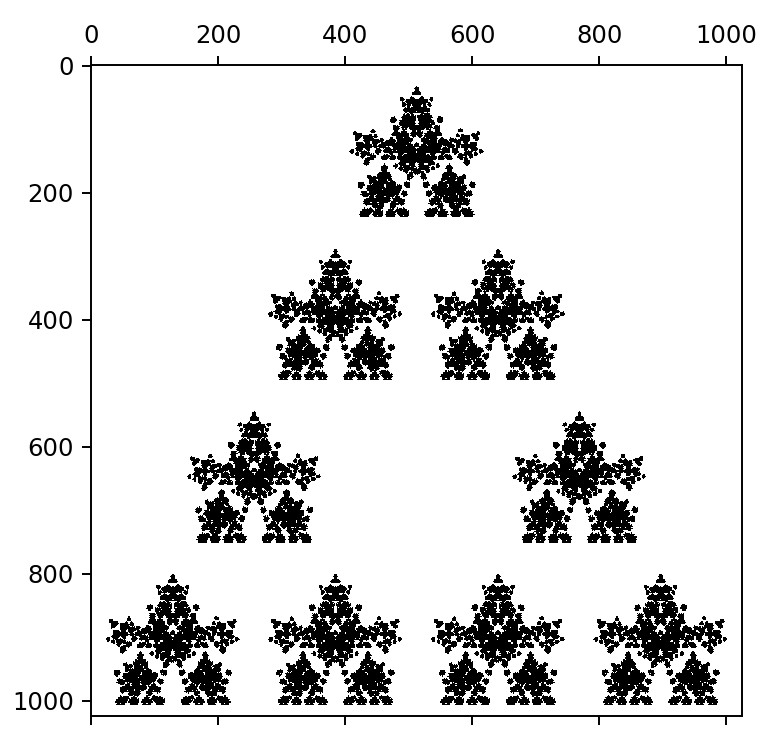} \quad 
	 	\includegraphics[width=4cm]{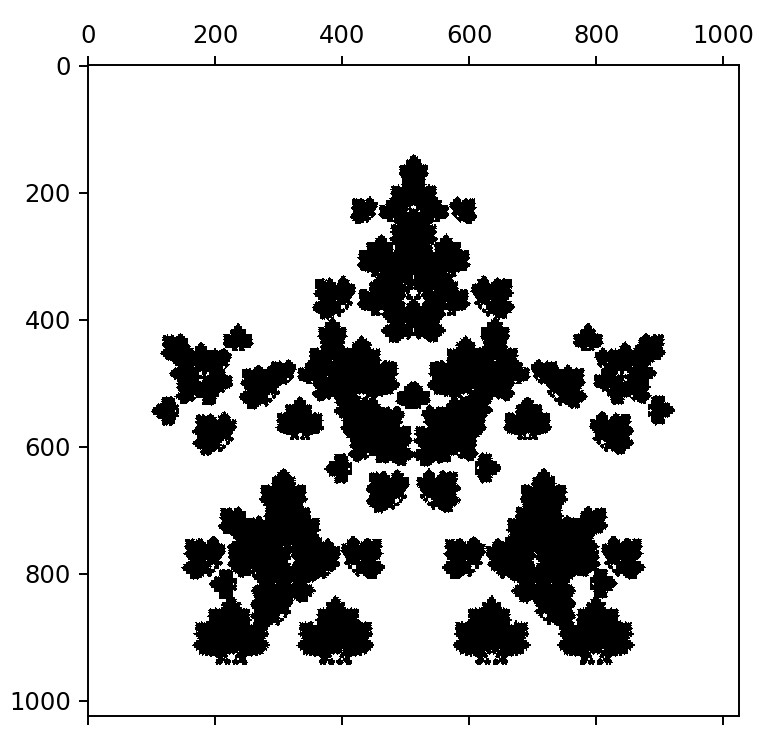} \quad 
	 	\includegraphics[width=4cm]{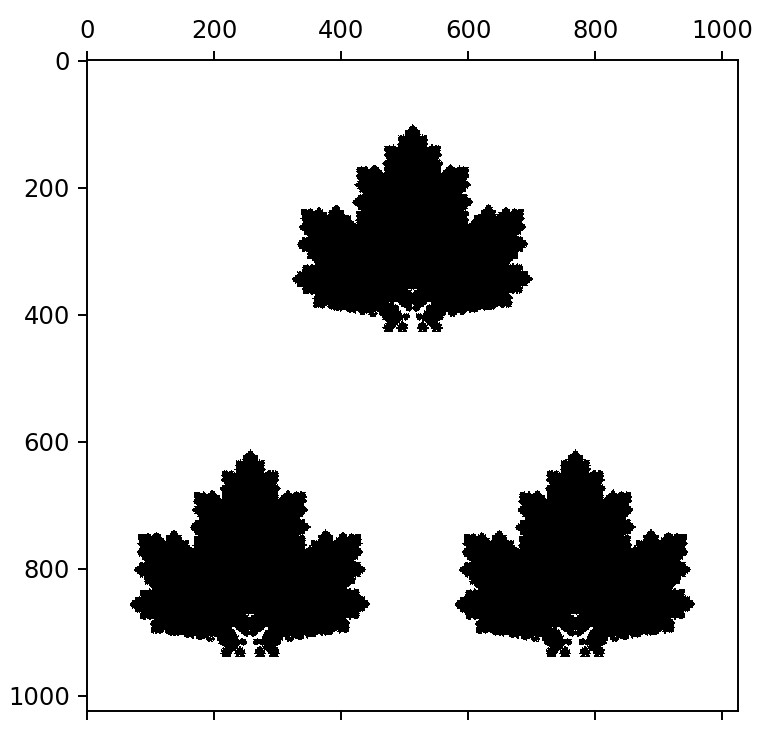}
	 	\caption{From left to the right, Sierpi\'nski and Maple blends: (1) , (2) and (3).}
	 	\label{fig:blend sierp mapl}
	 \end{figure}
	 Blends from Figure~\ref{fig:blend sierp mapl}:
	 \begin{enumerate}
	 	\item $\theta=(\boxed{1, 1}, 2, 1, 2, 1, 1, 2, 2, 1, 1, 2, 2, 2, 2, 1, 1, 1, 1, 1)$,  $\beta(\theta,1) \approx 1.3163$ and $\beta(\theta,2) \approx 1.9208$.
	 	\item $\theta=(2, \boxed{1}, 2, \boxed{1}, 2, 2, 2, 2, 2, 1, 1, 2, 2, 1, 2, 2, 1, 2, 2, 2)$,  $\beta(\theta,1) \approx 2.5778$ and $\beta(\theta,2) \approx 1.6048$.
	 	\item  $\theta=(1, \boxed{2, 2, 2, 2, 2, 2, 2}, 1, 2, 1, 2, 2, 1, 1, 1, 1, 1, 1, 1)$,  $\beta(\theta,1) \approx 2.6527116288$ and $\beta(\theta,2) \approx 1.5867172352$.
	 \end{enumerate}
	  We notice that in the blend (1) one have $\beta(\theta,1) \approx 1.3163 < \beta(\theta,2) \approx 1.9208$ which indicates that the attractor is much more Sierpi\'nski like that Maple like. On the other hand, in the blend (3) one have $\beta(\theta,1) \approx 2.6527116288> \beta(\theta,2) \approx 1.5867172352$ which indicates that the attractor is much more Maple like. Just in the last iteration one get  $\theta_1=1$ (meaning that we iterate last by the Sierpi\'nski operator $F_{R_{1}}$) after a series of 2's (meaning that we iterate last by the Maple operator $F_{R_{2}}$), but it only creates 3 copies of a Maple like miniature.  Finally, in the blend (2) most of the last iterations are by $F_{R_{2}}$, including the last one, except by $\theta_2=\theta_4=1$, so the general form is a deformed maple.
\end{example}

\begin{example} \label{examp: sierp maple and p2}
	  Consider $X=[0,1]^2$ and the IFSs $R_1=(X, f_1, f_2,f_3)$ and  $R_2=(X, f_1, f_2,f_3,f_4)$ from Example~\ref{examp: sierp and maple} , and in addition we introduce $R_3=(X, f_1, f_2,f_3,f_4)$ given by
	  \[\begin{cases}
	  	f_1(x,y)= \left( \tfrac{1}{3}x + \tfrac{1}{4}y, \tfrac{1}{12}x + \tfrac{19}{48}y \right)\\
	  	f_2(x,y)= \left( \tfrac{1}{3}x + \tfrac{1}{4}y, -\tfrac{1}{12}x + \tfrac{13}{48}y + \tfrac{1}{2} \right)\\
	  	f_3(x,y)=  \left( \tfrac{1}{3}x - \tfrac{1}{4}y + \tfrac{1}{2}, \tfrac{1}{12}x + \tfrac{13}{48}y + \tfrac{1}{8} \right)\\
	  	f_4(x,y)=\left( \tfrac{1}{3}x - \tfrac{1}{4}y + \tfrac{1}{2}, -\tfrac{1}{12}x + \tfrac{19}{48}y + \tfrac{3}{8} \right)
	  \end{cases}
	  \]
	  As before $A_{R_1}$ is a Sierpi\'nski like attractor (see Figure~\ref{fig:sierp}),   $A_{R_2}$ the Maple leaf attractor (see Figure~\ref{fig:sierp}) and $A_{R_3}$  is the attractor depicted in  Figure~\ref{fig:new} (see \cite{Oli17} for more information on this IFS, which appear as the power two of a non-contractive IFS in connection with GIFS attractors, see also \cite{MicMih08}). We consider the blend of these three attractors. It is well-known that $\lambda_{R_1}=0.5$,  $\lambda_{R_2}=0.8$  and $\lambda_{R_3}=0.5435$  thus we still have  $\lambda_{\mathcal{R}}=0.8$.   
	  \begin{figure}[H]
	  	\centering
	  	\includegraphics[width=6cm]{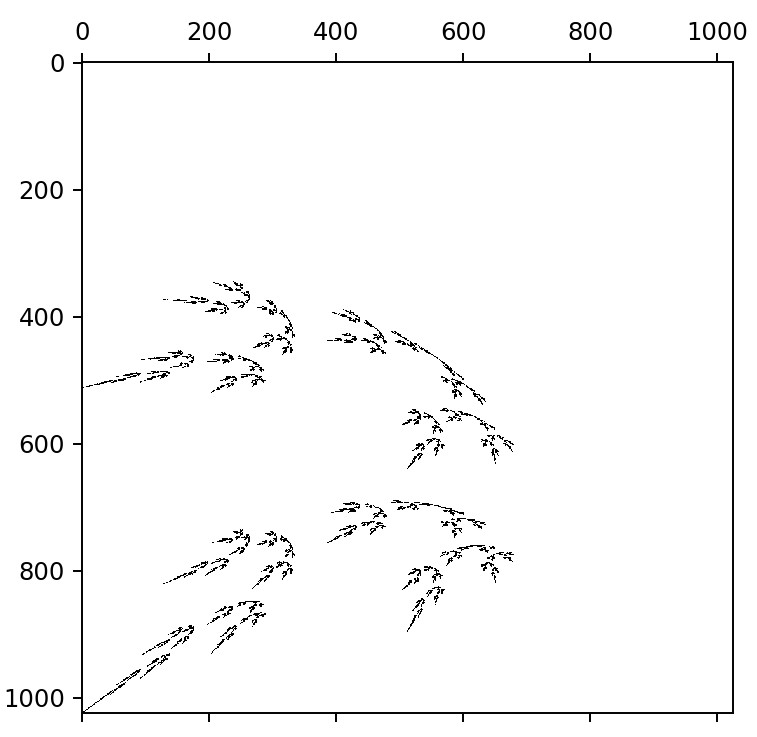} 
	  	\caption{Attractor $R_3$.}
	  	\label{fig:new}
	  \end{figure}	  
	  In order to apply the Algorithm~\ref{fig:blend approx algorithm} we consider the same parameters as in Example~\ref{examp: sierp and maple}. As $\lambda_{\mathcal{R}}=0.8$ remains the same, the error $d_H(\mathcal{A}(\theta), Y_1)$ is still $\approx 0.0198$ at most. 
	  
	  We now exhibit some blend discretizations $Y_1$, together with the blend sequence that generate it, and the blend coefficients. 
	  \begin{figure}[H]
	  	\centering
	  	\includegraphics[width=5cm]{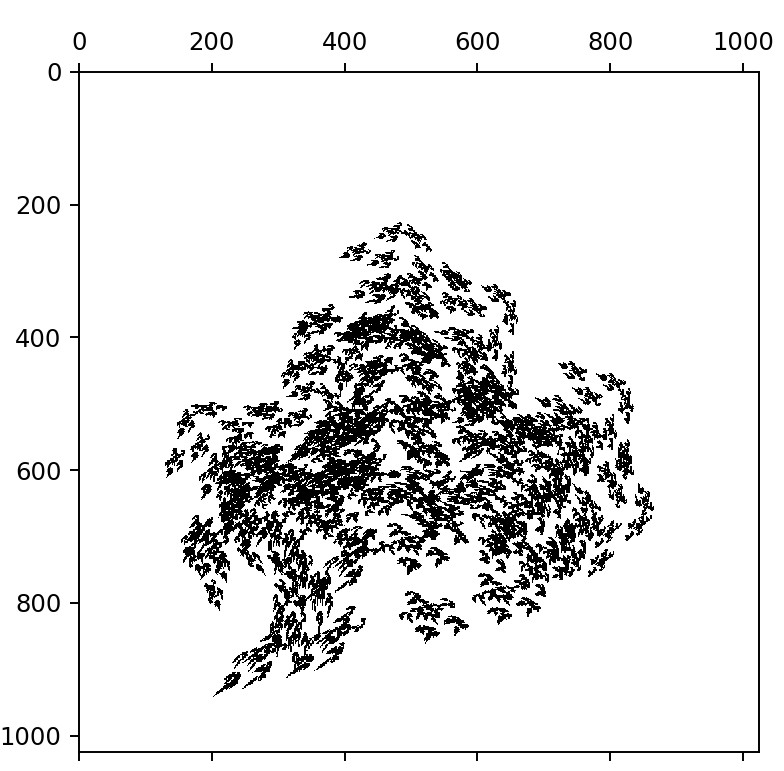} \quad 
	  	\includegraphics[width=5cm]{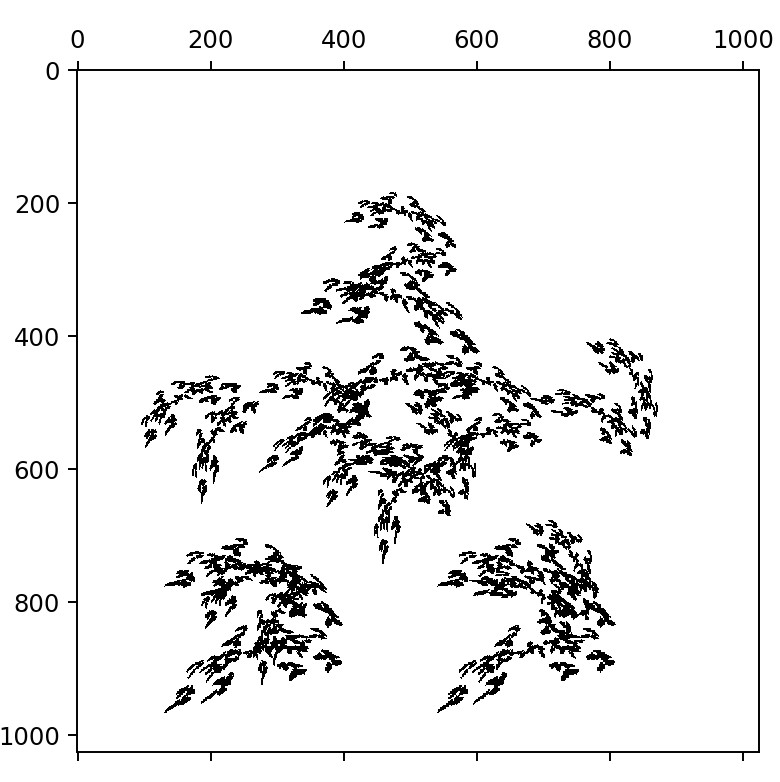} \\
	  	\includegraphics[width=5cm]{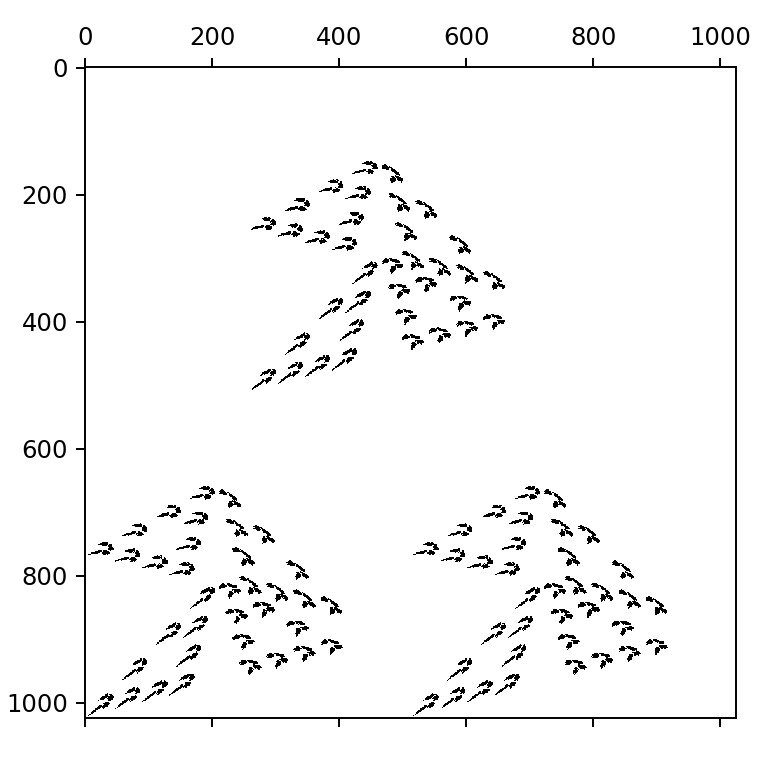}\quad
	  	\includegraphics[width=5cm]{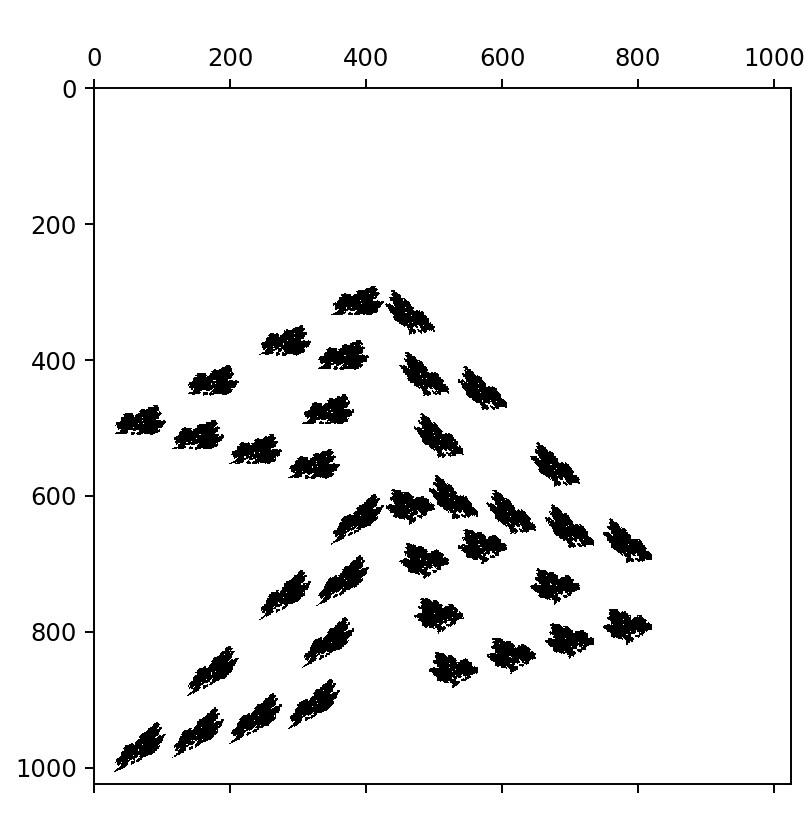} \\
	  	\includegraphics[width=5cm]{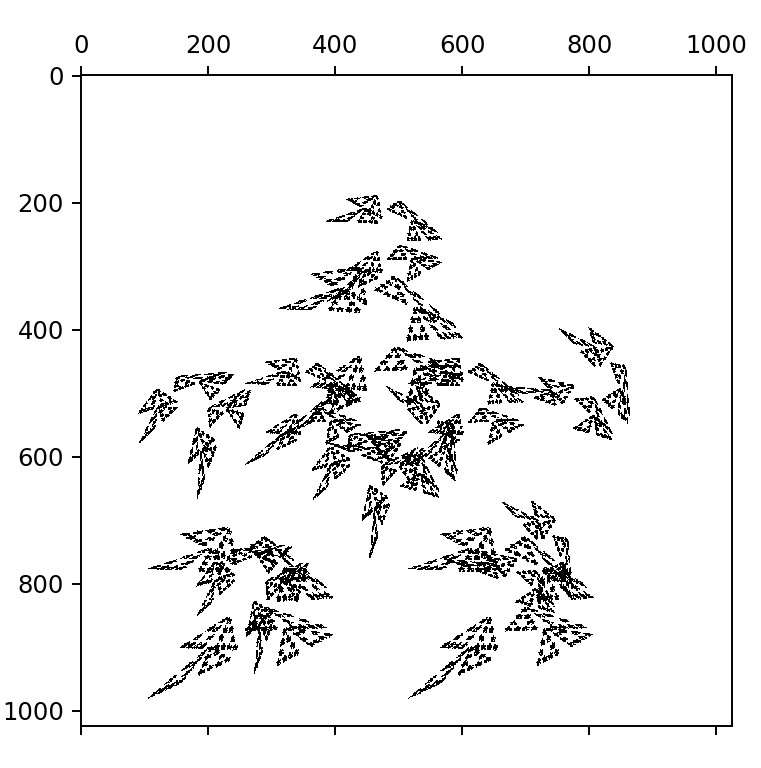} \quad 
	  	\includegraphics[width=5cm]{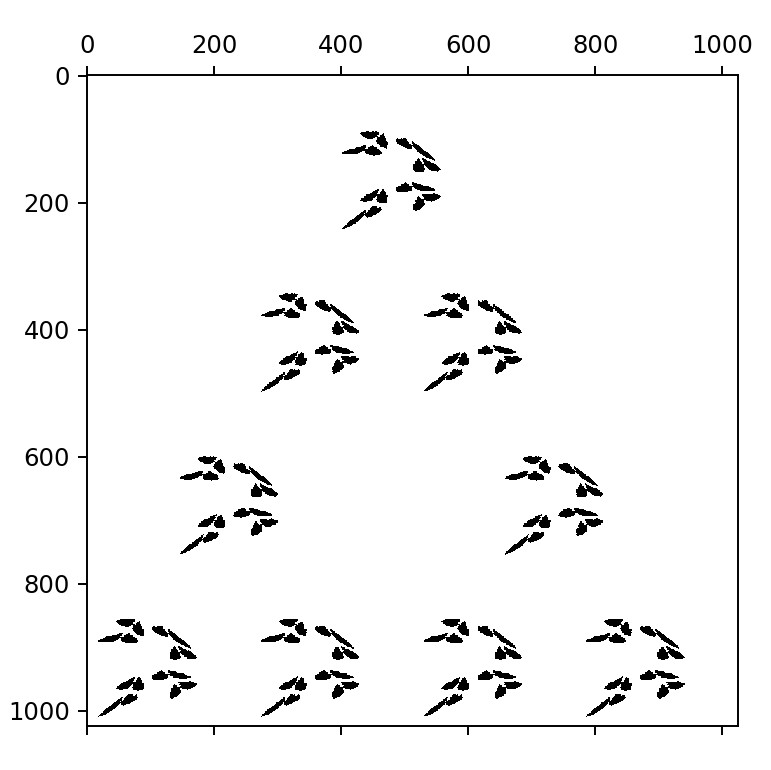}\\
	  	\includegraphics[width=5cm]{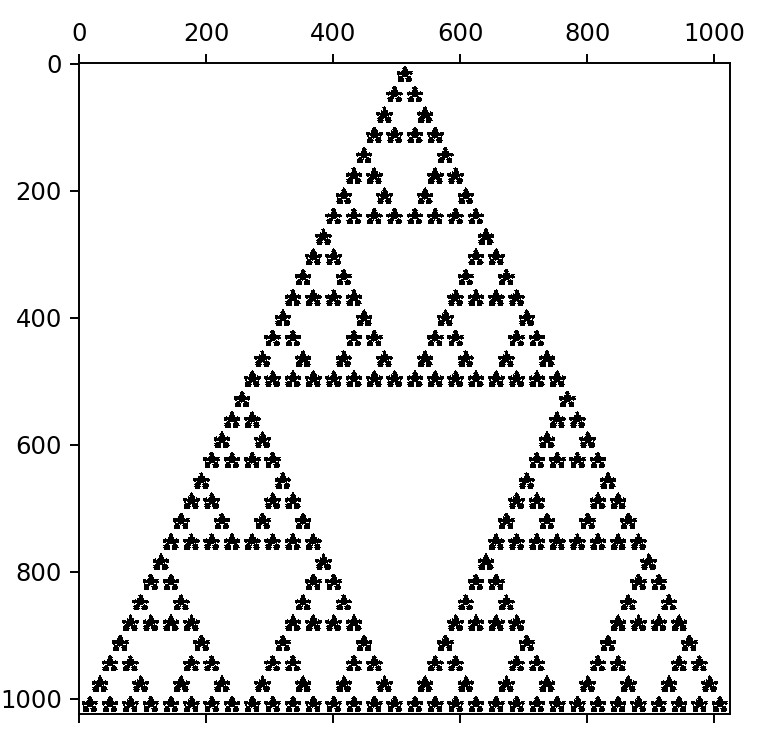} \quad 
	  	\includegraphics[width=5cm]{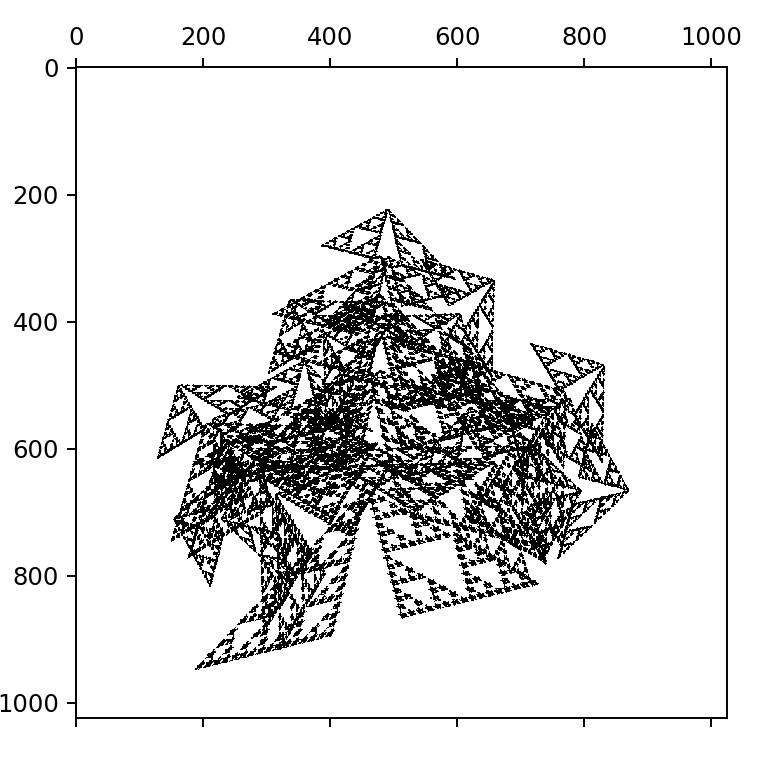} \\
	  	\caption{From left to the right, and up to bottom,  blends (1) -- (8).}
	  	\label{fig:blend sierp mapl p2}
	  \end{figure}
	  Blends from Figure~\ref{fig:blend sierp mapl p2}:
	  \begin{enumerate}
	  	\item $\theta=(2, 2, 3, 1, 2, 1, 3, 2, 3, 1, 1, 2, 1, 3, 3, 1, 2, 2, 1, 2)$,  $\beta(\theta,1) \approx 3.0165$, $\beta(\theta,2) \approx 1.6612$  and $\beta(\theta,3) \approx 2.8708$.
	  	\item $\theta=(2, 1, 3, 2, 1, 3, 2, 2, 2, 2, 3, 3, 2, 1, 2, 1, 3, 2, 2, 2)$,  $\beta(\theta,1) \approx 2.3743$, $\beta(\theta,2) \approx 1.7715$  and $\beta(\theta,3) \approx 2.5830$. 
	  	
	  	\item $\theta=(1, 3, 1, 1, 3, 2, 3, 2, 3, 2, 2, 3, 2, 2, 3, 2, 2, 2, 2, 3)$,  $\beta(\theta,1) \approx 1.3930$, $\beta(\theta,2) \approx 2.0389$  and $\beta(\theta,3) \approx 1.7617$.
	  \item $\theta=(3, 1, 1, 2, 2, 2, 3, 1, 2, 3, 1, 3, 1, 1, 2, 2, 2, 1, 2, 2)$,  $\beta(\theta,1) \approx 1.8734$, $\beta(\theta,2) \approx 2.0242$  and $\beta(\theta,3) \approx 1.7141$. 
	  
	  \item $\theta=(2, 1, 3, 3, 1, 1, 2, 1, 3, 2, 3, 1, 2, 2, 1, 3, 2, 2, 3, 1)$,  $\beta(\theta,1) \approx 2.1762$, $\beta(\theta,2) \approx 1.8474$  and $\beta(\theta,3) \approx 2.3334$.  
	  \item $\theta=(1, 1, 3, 3, 2, 2, 2, 2, 2, 2, 2, 2, 2, 2, 2, 2, 2, 2, 2, 1)$,  $\beta(\theta,1) \approx 1.4947$, $\beta(\theta,2) \approx 1.9610$  and $\beta(\theta,3) \approx 2.0362$.

	  \item $\theta=(\boxed{1, 1, 1, 1, 1, 2, 1}, 2, 2, 1, 1, 3, 2, 3, 3, 2, 3, 1, 3, 3)$,  $\beta(\theta,1) \approx 1.0460$, $\beta(\theta,2) \approx 1.9892$  and $\beta(\theta,3) \approx 2.0313$. 
	  \item $\theta=(2, 2, 3, 1, 1, 1, 1, 2, 1, 1, 3, 3, 1, 2, 2, 1, 3, 1, 1, 1)$,  $\beta(\theta,1) \approx 2.8099$, $\beta(\theta,2) \approx 1.6916$  and $\beta(\theta,3) \approx 2.7984$. 

	  \end{enumerate}
	  We notice that in the blend (7) one have  $\beta(\theta,1) \approx 1.0460$, very close to the minimum, which indicates that the attractor is much more Sierpi\'nski  like. Indeed, most of the last iterations are by $F_{R_{1}}$, including the last one, except by $\theta_6=2$.  However, for more than two IFSs the influence of each one in the blend is much more elusive to decipher. 
\end{example}

\subsection{Application: Canright's Envelope}

From the original work of \cite{Can94}, later improved by \cite{AnMi23}, one can find a  bound for the attractor $\mathcal{A}_{\mathcal{R}}$ of an IFS $\mathcal{R}$ by a set of balls centered at the fixed points of each map in $\mathcal{R}$, called Canright's Envelope of $\mathcal{A}_{\mathcal{R}}$.  In the recent work of $\mathcal{A}_{\mathcal{R}}$ the improved the formula for the radii of these balls and developed an extension allowing to find finer coverings, but with a larger number of balls appearing from the fixed points of powers of the original IFS. As an application to the blend of IFS attractors, we can use the covering to estimate the maximum distance from each blend to the individual attractors.

From Definition~\ref{def: blending IFS},  if $(X,d)$ is a complete metric space and $R_{i}=(X, f_{j}^{i}, j=1,\ldots,n_{i} )$, for $i=1,\ldots,N$ a family of contractive IFSs one consider the IFS $\mathcal{R}:=( K^*(X), F_{R_{i}}, i=1,\ldots,N )$ as  the blending of the IFSs $R_{i}$. 

From \cite{Can94}, if $A_{R_{i}}, i=1,\ldots,N$ are the attractors of the respective $R_{i}$, that is, $F_{R_{i}}(A_{R_{i}})=A_{R_{i}}$, and $t_i$ are the solutions of
\begin{equation}\label{eq:canright }
	t_i=\lambda_{R_{i}} \left( \max_{j \neq i} d_H(A_{R_{i}}, A_{R_{j}})   + t_j\right) 
\end{equation}
then $\mathcal{A}_{\mathcal{R}} \subset \bigcup_{j=1,\ldots,N} B(A_{R_{j}}, t_j):=E$, where  \[B(A_{R_{j}}, t_j):=\{Q \in K^*(X) \;|\; d_H(Q, A_{R_{j}})<t_j\}.\] The set $E$ is called   Canright's Envelope of $\mathcal{A}_{\mathcal{R}}$. As pointed in \cite{AnMi23} there is no practical way to compute these numbers.

However,  \cite[Proposition 3.3]{AnMi23} provides a new formulation, where the numbers $r_i$ are the solutions of
\begin{equation}\label{eq:canright ang mic}
	r_i=\lambda_{R_{i}} \left( \max_{j , i} d_H(A_{R_{i}}, A_{R_{j}})   +  \max_{j \neq i} r_j\right) 
\end{equation}
satisfy $\mathcal{A}_{\mathcal{R}} \subset \bigcup_{j=1,\ldots,N} B[A_{R_{j}}, r_j]:=\mathcal{C}_1$, where  \[B[A_{R_{j}}, t_j]:=\{Q \in K^*(X) \;|\; d_H(Q, A_{R_{j}})\leq t_j\}.\]  The set $\mathcal{C}_1$ is called (degree 1) covering of $\mathcal{A}_{\mathcal{R}}$ and provide an explicit solution of Equation~\eqref{eq:canright ang mic} in the ordered form $r_{i_{1}} \leq r_{i_{2}} \leq \ldots r_{i_{N}} $:
\begin{figure}[H]
	\centering
	\includegraphics[width=8cm]{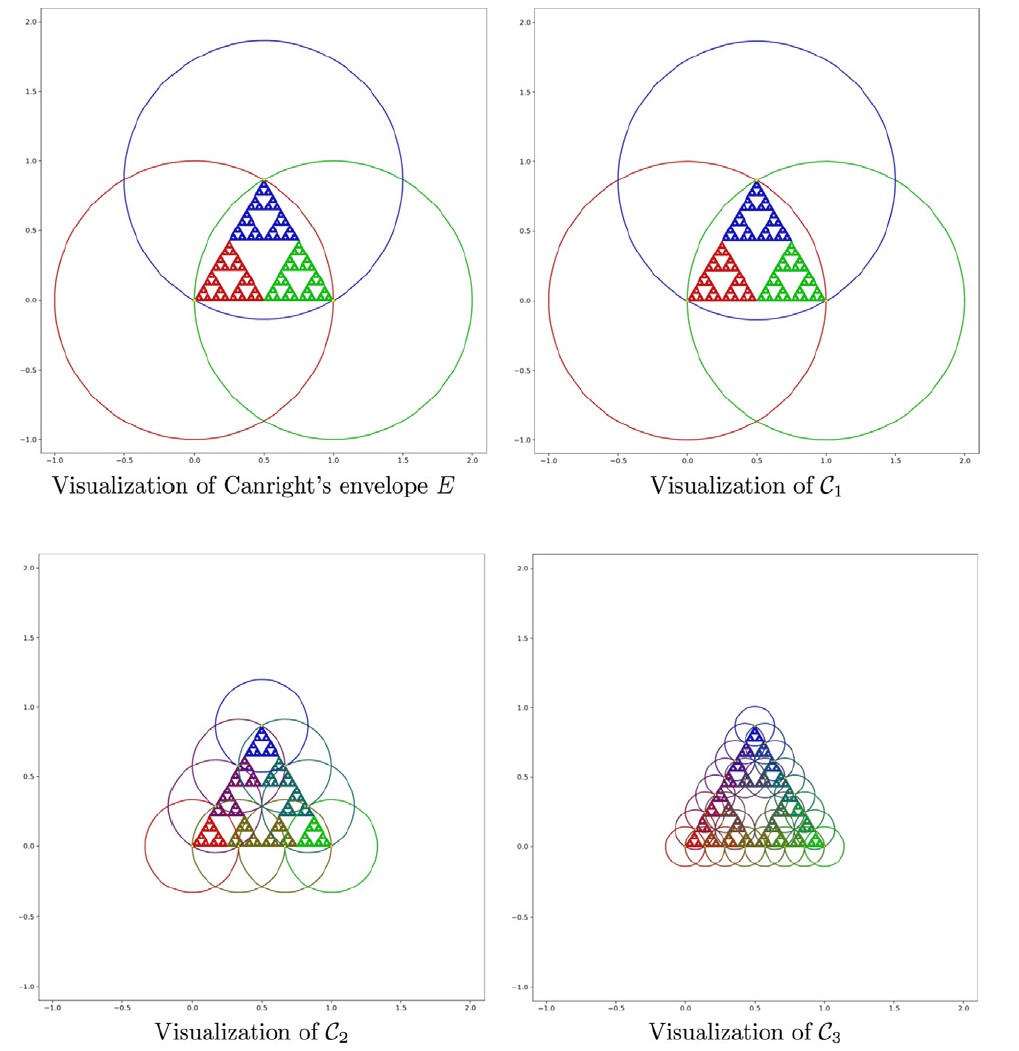} 
	\caption{Attractor $R_1$(Sierpi\'nski’s Triangle)  coverings from \cite[Example 4.2]{AnMi23}.}
	\label{fig:coverings}
\end{figure}	  
\begin{theorem}\cite[ Theorem 3.1]{AnMi23} \label{thm:3.1 Ang and Micu} 
	In the above-mentioned framework, we have
	\[r_{i_{j}}:= M \lambda_{R_{i_j}} \frac{1+\lambda_{R_{i_N}}}{1-\lambda_{R_{i_{N-1}}} \lambda_{R_{i_N}}}, \; j=1,\ldots,N-1\]
	and 
	\[r_{i_{N}}:= M \lambda_{R_{i_N}} \frac{1+\lambda_{R_{i_{N-1}}}}{1-\lambda_{R_{i_{N-1}}} \lambda_{R_{i_N}}},\]
	where $M:=\max_{j , i} d_H(A_{R_{i}}, A_{R_{j}})$.
\end{theorem}

\begin{example}\label{examp: applic Canrights}   We consider the IFS $R_1, R_2$ and $R_3$ given in Example~\ref{examp: sierp maple and p2} where $N=3$ and $\lambda_{R_1}=0.5$,  $\lambda_{R_2}=0.8$  and $\lambda_{R_2}=0.5435$.
	
	Using one of several methods to evaluate Hausdorff distance one can compute 
	\[M=\max(d_H(A_{R_{1}}, A_{R_{2}}) , \; d_H(A_{R_{1}}, A_{R_{3}}), \; d_H(A_{R_{2}}, A_{R_{3}})),\]
	and the Equation~\ref{eq:canright ang mic} became
	\[\begin{cases}
		r_1=\lambda_{R_{1}} \left( M  +   \tilde r \right)=0.5\, \left( M  +   \tilde r \right)\\
		r_2=\lambda_{R_{2}} \left( M  +   \tilde r \right)=0.8\, \left( M  +   \tilde r \right)\\
		r_3=\lambda_{R_{3}} \left( M  +   \tilde r \right)=0.5435\, \left( M  +   \tilde r \right)
	\end{cases} 	
	\] 
	Since $\tilde{r}=\max(r_1,r_2,r_3) $ we obtain $i_1=1, i_2=3\text{ and }i_3=3$. Hence,
	\[\tilde{r}=\max(r_1,r_2,r_3)  = 0.8\, \left( M  +   \tilde r \right)\]
	or equivalently
	\[\tilde{r}=M \frac{0.8}{1-0.8}= 4 M.\]
	Consequently, $r_1= 0.5\, \left( M  +  4 M \right)= \frac{5}{2} \, M=2.5 \,M$, $r_2=\tilde{r}= 4\, M$ and $r_3=0.5435\, \left( M  +   4 M \right)=2.7175\, M$. 
	
	In conclusion, all blend $\mathcal{A}(\theta)$ is, at most $2.5 \,M$ from the Sierpi\'nski $A_{R_{1}}$,   at most $4 \,M$ from the Maple leaf $A_{R_{2}}$ or  at most $2.7175 \,M$ from the attractor $A_{R_{3}}$,  with respect to the Hausdorff distance.
	
	A free approximation (in some references $d_H(\text{Maple leaf, Sierpi\'nski})\approx 0.289971$ can be founded, but it depends on the particular choice of functions, and we are not interested in the optimal computation, neither in provide error estimates for it)  provides 
	\[M=\max(d_H(A_{R_{1}}, A_{R_{2}}) , \; d_H(A_{R_{1}}, A_{R_{3}}), \; d_H(A_{R_{2}}, A_{R_{3}}))=\]
	\[=	\max(0.3123 , \; 0.4101, \; 0.3102) = 0.41+.\]
	Assuming this estimation one obtain, for any $\theta \in \Omega$, that
	$d_H(\mathcal{A}(\theta), A_{R_{1}}) \leq 2.5 \,M\approx 1.02$ or $ d_H(\mathcal{A}(\theta), A_{R_{2}}) \leq 4 \,M\approx  1.64 $ or $ d_H(\mathcal{A}(\theta), A_{R_{3}}) \leq 2.7175 \,M \approx  1.11.$ 	None of them is particularly good since $\operatorname{diam}(X)=\sqrt{2}\approx 1.41$.
	
	We could improve this bounds using, for instance, a degree 2 covering from \cite{AnMi23} (see Figure~\ref{fig:coverings}), but it will be required a great computational power to deal with nine fixed points combining the original attractors, to be precise, {\small  $\mathcal{A}(1,1,1,1,\ldots),$ $ \mathcal{A}(2,2,2,2,\ldots),$ $  \mathcal{A}(3,3,3,3,\ldots), $ $ \mathcal{A}(1,2,1,2,\ldots), $ $ \mathcal{A}(2,1,2,1,\ldots), $ $ \mathcal{A}(2,3,2,3,\ldots),$ $  \mathcal{A}(1,3,1,3,\ldots), $ $ \mathcal{A}(3,1,3,1,\ldots) $  and $\mathcal{A}(3,2,3,2,\ldots)$}. Note that the first three are just 
	$A_{R_{1}}$, $A_{R_{2}}$ and $A_{R_{3}}$.
\end{example}

\subsection{Application: Blending  invariant probabilities}
An immediate unfolding of our approach is the consideration of IFSp, which stands for IFSs with probabilities.  If  $(X,d)$ is a compact metric space and $R_{i}=(X, f_{j}^{i}, p_j^i, j=1,\ldots,n_{i} )$, for $i=1,\ldots,N$ a family of contractive IFSs with probabilities $\sum_{j=1,\ldots,n_{i} } p_j^i  =1$, one consider the IFS $\mathcal{S}:=( \mathcal{P}(X), M_{R_{i}}, i=1,\ldots,N )$ where $\mathcal{P}(X)$ is the set of Borel probabilities over $X$, endowed with the Monge-Kantorovich distance, and $M_{R_{i}}: \mathcal{P}(X) \to \mathcal{P}(X)$ is the Markov operator, associated to $R_{i}=(X, f_{j}^{i}, p_j^i, j=1,\ldots,n_{i} )$, given by 
\[M_{R_{i}}(\mu)(B):= \sum_{j=1,\ldots,n_{i} } p_j^i \mu((f_{j}^{i})^{-1}(B)),\]
for any Borelian $B \subseteq X$.

Under this simple conditions (see \cite{Hut81} or \cite{BarnDem85}),  each map $M_{R_{i}}$ is a contraction and has a unique fixed point $\mu_{R_{i}} \in \mathcal{P}(X)$, denoted the Hutchinson measure, or invariant measure, of the IFS $R_{i}$.

On the other hand $\mathcal{S}:=( \mathcal{P}(X), M_{R_{i}}, i=1,\ldots,N )$ will have an attractor, $A_{\mathcal{S}} \subset \mathcal{P}(X)$ and for each  blending sequence $\theta \in  \Omega$ the probability
\[\mu(\theta) :=\pi_{\mathcal{S}}(\theta) = \lim_{k \to \infty}  M_{R_{\theta_1}}(\cdots( M_{R_{\theta_k}}(\nu))  \in \mathcal{P}(X)\]
for a fixed $\nu \in \mathcal{P}(X)$, given by Theorem~\ref{thm: code map level 1}, is called the blend  of the invariant measures $\mu_{R_{i}}, i=1,\ldots,N$ by $\theta$.

It is worth to mention that, in this formulation, $\mathcal{S}$ is a deterministic IFS, despite the fact that each IFSp is probabilistic. Introduce probabilities in $\mathcal{S}$ would generate a second level of randomness since our blending is always deterministic.

Many questions arrive regarding the statistical properties of $\mu(\theta)$ in comparison with each $\mu_{R_{i}}, i=1,\ldots,N$, which are known to be ergodic even for place dependent probabilities cf.~\cite{Elton1987}. It should the subject of further research.

\end{document}